\documentclass[10pt,twoside]{article}
\usepackage{amssymb,amsmath,amsthm,enumerate,geometry}
\usepackage[mathscr]{eucal}
\input amssym.def
\newtheorem{example}{Example}[section]
\newtheorem{theorem}{Theorem}[section]
\newtheorem{lemma}{Lemma}[section]
\newtheorem{proposition}{Proposition}[section]
\newtheorem{remark}{Remark}[section]
\newtheorem{corollary}{Corollary}[section]
\newtheorem{definition}{Definition}[section]

\def\bb{\Bbb}
\def\bpr{\begin{proof}}
\def\epr{\end{proof}}
\def\be{\begin{equation}}
\def\ee{\end{equation}}
\def\bea{\begin{eqnarray}}
\def\eea{\end{eqnarray}}
\def\bean{\begin{eqnarray*}}
\def\eean{\end{eqnarray*}}
\def\ba{\begin{abstract}}
\def\ea{\end{abstract}}
\def\bt{\begin{theorem}}
\def\et{\end{theorem}}
\def\bl{\begin{lemma}}
\def\el{\end{lemma}}
\def\bp{\begin{proposition}}
\def\ep{\end{proposition}}
\def\br{\begin{remark}}
\def\er{\end{remark}}
\def\bc{\begin{corollary}}
\def\ec{\end{corollary}}
\def\bd{\begin{definition}}
\def\ed{\end{definition}}
\def\non{\nonumber}

\def\bs{\bigskip}

\begin{document}
 \baselineskip 20pt
\title{ Multi-window dilation-and-modulation frames on the half real line\thanks{Supported by the National Natural Science
Foundation of China (Grant No. 11271037).
\newline \hspace{1cm}}}
\author{Yun-Zhang Li$^{1}$\,\,\,\,\ \ Wei Zhang$^{2}$\\
College of Applied Sciences,\\ Beijing University of Technology,
Beijing 100124,  P. R. China\\E-mail:  1. yzlee@bjut.edu.cn\\2. zhangweihappy@emails.bjut.edu.cn}\date{}
\maketitle{}
\begin{abstract} Wavelet and Gabor systems are based on translation-and-dilation and translation-and-modulation operators, respectively. They have been extensively studied. However, dilation-and-modulation systems have not, and they cannot be derived from wavelet or Gabor systems. In this paper, we investigate a class of dilation-and-modulation systems in the causal signal space $L^{2}(\bb R_{+})$. $L^{2}(\bb R_{+})$ can be identified  a subspace of $L^{2}(\bb R)$ consisting of all $L^{2}(\bb R)$-functions supported  on $\bb R_{+}$, and is unclosed under the Fourier  transform. So the Fourier transform method
 does not work in $L^{2}(\bb R_{+})$. In this paper, we introduce the notion of $\Theta_{a}$-transform in $L^{2}(\bb R_{+})$, using $\Theta_{a}$-transform we characterize dilation-and-modulation frames and dual frames in $L^{2}(\bb R_{+})$; and present an explicit expression of all duals with the same structure for a general dilation-and-modulation frame for $L^{2}(\bb R_{+})$.  Interestingly, we prove that  an arbitrary frame of this form is always nonredundant whenever the number of the generators is $1$, and is always redundant whenever it is greater than $1$. Some examples are also provided to illustrate the generality of our results.
 \end{abstract}

\noindent {\bf 2010 MS Classification}:  42C40, 42C15

\section{Introduction}
\setcounter{equation}{0}

It is well known that  translation, modulation and  dilation are fundamental operations in wavelet analysis. The  translation
operator  $T_{x_{0}}$, the modulation operator $M_{x_{0}}$ with $x_{0}\in \bb R$ and the dilation operator $D_{c}$ with $0<c\ne 1$ are  defined by
$$T_{x_{0}}f(\cdot)=f(\cdot-x_{0}),{\mbox{ }}
M_{x_{0}}f(\cdot)=e^{2\pi ix_{0}\cdot}f(\cdot){\mbox{ and }}
D_{c}f(\cdot)=\sqrt{c}f(c\cdot)$$ for $f\in L^{2}(\bb R)$,
respectively. Given  a finite subset $\Psi$ of $L^{2}(\bb R)$,
Gabor frames of the form
 \be\label{gab}\{M_{mb}T_{na}\psi:\,m,\,n\in \bb Z,\,\psi\in \Psi\}\ee
 and  wavelet frames of the form
\be\label{wav}\{D_{a^{j}}T_{bk}\psi:\,j,\,k\in \bb Z, \,\psi\in \Psi\}\ee
  with $a$, $b>0$  have been extensively studied ([\ref{O}, \ref{FS1}, \ref{FS2}, \ref{Gro01}, \ref{hei11}-\ref{her96}]).
However, dilation-and-modulation frames of the form
\be\label{md}
\{M_{mb}D_{a^{j}}\psi:\,m,\,j\in \bb Z,\,\psi\in \Psi\}{\mbox{ with }}a,\mbox{ }b>0
\ee
have not. Observe that the Fourier transform of (\ref{md}) is
\be\label{fouri}\{T_{mb}D_{a^{j}}\hat\psi:\,m,\,j\in \bb Z,\,\psi\in \Psi\}\ee
It does not fall into the framework of the above wavelet and Gabor systems.
Our focus in this paper will be on a class of dilation-and-modulation frames
for $L^2(\bb R_+)$ with  $\bb R_+=(0,\,\infty)$.
$L^2(\bb R_+)$ can be considered as a closed subspace of $L^2(\bb R)$ consisting of all
functions in $L^2(\bb R)$ which vanish outside $\bb R_+$. And it  can model causal signal space. In practice, time variable cannot be negative.

For subspace Gabor and wavelet frames of the forms (\ref{gab}) and (\ref{wav}) respectively, we refer to [\ref{CC01}, \ref{DDG}-\ref{D1}, \ref{GH01}-\ref{GL}, \ref{Gh09}, \ref{Gh11}, \ref{jiali14}, \ref{jiali15},
\ref{li15}, \ref{li152}-\ref{lizhou11},  \ref{s93}, \ref{volk95}, \ref{zh14}, \ref{lizhou10}] and references therein for details.  It is easy to check that
there exists no nonzero function $\psi$ such that
$$T_{nc}\psi(\cdot)=0{\mbox{ on }}(-\infty,\,0)$$
for some $c>0$ and all $n\in\bb Z$. This implies that $L^2(\bb R_+)$ admits
no frame of the form (\ref{gab}), (\ref{wav}) or (\ref{fouri}).  So it is natural to ask
how we construct frames for $L^2(\bb R_+)$ with good structures.
Two methods are known to us for this purpose. One  is to construct frames for $L^2(\bb R_+)$  consisting of  a subsystem of (\ref{wav}) and some inhomogeneous refinable function-based ``boundary wavelets". For details, we refer to [\ref{chen99}, \ref{dah00},  \ref{jia00},  \ref{jia01},  \ref{li00}, \ref{mic94},  \ref{str98},  \ref{sun99}] and references therein. The other is
to use the Cantor  group operation  and Walsh series theory to introduce the
notion of (frame) multiresolution analysis in $L^2(\bb R_+)$, and then derive wavelet frames similarly to the case of $L^2(\bb R)$. For details, we refer to [\ref{al10}, \ref{ev15}-\ref{far11},  \ref{lang96},  \ref{lang98}, \ref{ni04},  \ref{ni00},  \ref{sh12}] and references therein.
The references [\ref{gms98}] and [\ref{gms01}] also have something to do with
this problem. In [\ref{gms98}],
numerical experiments were made  to establish that the nonnegative integer shifts of the
Gaussian function form a Riesz sequence in $L^{2}(\bb R_{+})$. And in [\ref{gms01}], a sufficient condition was obtained to determine
whether or not the nonnegative translates  of a given function form a Riesz sequence on $L^{2}(\bb R_{+})$.

Given $a>1$, a measurable function $h$ defined on $\bb R_+$ is said to be $a$-{\it{dilation periodic}} if $h(a\cdot)=h(\cdot)$ a.e. on $\bb R_+$.  Throughout this paper, we  denote by $\{\Lambda_{m}\}_{m\in\bb Z}$ the sequence of
$a$-dilation periodic functions defined by
\be \label{1111}
\Lambda_{m}(\cdot)=\frac{1}{\sqrt{a-1}}e^{\frac{2\pi im\cdot}{a-1}}\mbox{ on }[1,\,a)\mbox{ for each }m\in\bb Z.
\ee
Motivated by the above works, we in this paper investigate  the dilation-and-modulation systems in $L^2(\bb R_+)$ of the form:
\be \label{000}
{\cal MD}(\Psi,\,a)=\left\{\,{\Lambda_{m}}D_{a^{j}}\psi_l:\,m,\,j\in \bb Z, \,1\le l\le L\,\right\}
\ee
under the following General setup:

\bs
{\bf General setup:}

(i) $a$ is a fixed positive number greater than 1.

(ii) $\Psi=\{\psi_1,\,\psi_2,\,\cdots,\,\psi_L\}$ is a finite subset of $L^2(\bb R_+)$ with cardinality $L$.

\bs

For $\Phi=\{\varphi_1,\,\varphi_2,\,\cdots,\,\varphi_L\}$, we define ${\cal MD}(\Phi,\,a)$ similarly to (\ref{000}). The system
${\cal MD}(\Psi,\,a)$ is slightly like   but differs from (\ref{md}).  The modulation factor
$e^{2\pi imb\cdot}$ in  (\ref{md}) is $\bb Z$-periodic according to addition, while $\Lambda_{m}$ in (\ref{000}) is $a$-dilation periodic.
 Let $a$ and $\Psi$ be as in the general setup. The system ${\cal MD}(\Psi,\,a)$ is called
a {\it frame} for $L^2(\bb R_+)$ if there exist
$0<C_1\le C_2<\infty$ such that \be\label{0000}C_1\|f\|_{L^{2}(\bb R_{+})}^2\le \sum\limits_{l=1}^{L}\sum\limits_{m,j\in \bb Z}\left|\langle f,\,\Lambda_{m}D_{a^{j}}\psi_l\rangle_{L^{2}(\bb R_{+})}\right|^{2}\le C_2\|f\|_{L^{2}(\bb R_{+})}^2
{\mbox{ for }}f\in L^2(\bb R_+),\ee  where $C_1$, $C_2$ are called {\it frame bounds}; it is called a
{\it Bessel sequence} in $L^2(\bb R_+)$
if the right-hand side inequality in (\ref{0000}) holds,
where $C_2$ is called a {\it Bessel bound}.
In particular, it is
called a {\it Parseval frame} if $C_1=C_2=1$ in
(\ref{0000}).  Given a frame ${\cal MD}(\Psi,\,a)$  for $L^2(\bb R_+)$,
a sequence ${\cal MD}(\Phi,\,a)$  is called  a {\it dual} (or an ${\cal MD}$-{\it
dual}) of  ${\cal MD}(\Psi,\,a)$ if it is a frame such that
 \be\label{ddual}f=\sum\limits_{l=1}^{L}\sum\limits_{m,j\in \bb Z}
 \langle f,\,\Lambda_{m}D_{a^{j}}\varphi_l\rangle_{L^{2}(\bb R_{+})}\Lambda_{m}D_{a^{j}}\psi_l{\mbox{ for  }}
f\in L^{2}(\bb R_{+}).\ee
It is easy to check that ${\cal MD}(\Psi,\,a)$ is also a
dual of ${\cal MD}(\Phi,\,a)$  if ${\cal MD}(\Phi,\,a)$
 is a dual of  ${\cal MD}(\Psi,\,a)$. So, in this case, we say ${\cal MD}(\Psi,\,a)$  and  ${\cal MD}(\Phi,\,a)$
form a pair of dual frames for  $L^{2}(\bb R_{+})$. By the knowledge of frame theory,
 ${\cal MD}(\Psi,\,a)$ and  ${\cal MD}(\Phi,\,a)$
form a pair of dual frames for  $L^{2}(\bb R_{+})$ if they are Bessel sequences and satisfy  (\ref{ddual}).  The fundamentals of frames can be found in [\ref{O},  \ref{duffin},   \ref{hei11}, \ref{RY}]. Observe that
$L^{2}(\bb R_{+})$ is the Fourier transform of the Hardy space $H^{2}(\bb R)$ which is a reducing subspace of $L^{2}(\bb R)$ defined by
$$H^{2}(\bb R)=\{f\in L^{2}(\bb R):\hat{f}(\cdot)=0{\mbox{ a.e. on}} (-\infty,\,0)\}.$$
 Wavelet frames in $H^{2}(\bb R)$ of the form  (\ref{wav})
   were studied in [\ref{her96}, \ref{s93}, \ref{volk95}]. By the Plancherel theorem, an $H^{2}(\bb R)$-frame
$\{D_{a^{j}}T_{bk}\psi:\,j,\,k\in \bb Z, \,\psi\in \Psi\}$ leads to  a $L^{2}(\bb R_{+})$-frame
\be \label{001}
\{e^{-2\pi ia^{j}k\cdot}\hat{\psi}(a^{j}\cdot):\,j,\,k\in \bb Z, \,\psi\in \Psi\,\}.
\ee
In (\ref{001}), $e^{-2\pi ia^{j}k\cdot}$ is $a^{-j}\bb Z$-periodic with respect to addition, and the period varies with  $j$. However,
$\Lambda_{m}$ in (\ref{000}) is $a$-dilation periodic, and unrelated to $j$. Therefore, the system (\ref{000}) differs from  (\ref{001}) for $L^{2}(\bb R_{+})$, and is of independent interest.  Actually, it is slightly related to a kind of  function-valued frames  in \cite{benli98}.

This paper focuses on the  theory of
$L^{2}(\bb R_{+})$-frames of the form (\ref{000}). It cannot be derived from
the well known wavelet and Gabor systems, and its operation is more intuitive
when compared with the Cantor  group and Walsh series-based systems in [\ref{al10}, \ref{ev15}-\ref{far11},  \ref{lang96},  \ref{lang98}, \ref{ni04},  \ref{ni00},  \ref{sh12}].   Also $L^{2}(\bb R_{+})$ is unclosed under the Fourier transform.
In particular,  the Fourier transform of a compactly supported nonzero function in $L^{2}(\bb R_{+})$ lies outside this space. Therefore, the Fourier
transform cannot be used in our setting, and we need to find a new method.

The rest of this paper is organized follows. In Section 2, we introduce the notion of $\Theta_{a}$-transform, and give a $\Theta_{a}$-transform domain
characterization of a dilation-and-modulation system ${\cal MD}(\Psi,\,a)$
being complete, a Bessel sequence and a frame in $L^{2}(\bb R_{+})$, respectively. In Section 3, using $\Theta_{a}$-transform we characterize dual frame pairs of the form $({\cal MD}(\Psi,\,a),\,{\cal MD}(\Phi,\,a))$, and obtain an explicit expression of all ${\cal MD}$-duals of a general frame ${\cal MD}(\Psi,\,a)$ for $L^{2}(\bb R_{+})$. We also prove that an arbitrary frame ${\cal MD}(\Psi,\,a)$ is a Riesz basis if and only if  $L=1$. It means that the frame ${\cal MD}(\Psi,\,a)$ is always nonredundant whenever $L=1$, and it is always redundant whenever $L>1$. In Section 4, we give some examples
of   ${\cal MD}$-dual frame pairs for $L^{2}(\bb R_{+})$ to illustrate the generality of our results. They show that the achieved results in this paper
 provide us with an easy method to construct  ${\cal MD}$-dual frames for $L^{2}(\bb R_{+})$ with the window functions having  good properties such as having bounded supports and certain smoothness.

\section{ $\Theta_{a}$-transform domain frame characterization}
\setcounter{equation}{0}

Let $a$ and $\Psi$ be as in the general setup. In this section, by introducing $\Theta_{a}$-transform we give the conditions of completeness, Bessel sequence and frame of ${\cal MD}(\Psi,\,a)$ in $L^{2}(\bb R_{+})$, respectively.

\bd\label{theta} Let $a$ be as in the general setup. For $f\in L^{2}(\bb R_{+})$, we define
\be \label{2001}
\Theta_{a}f(x,\xi)=\sum\limits_{l\in\bb Z}a^{\frac{l}{2}}f(a^{l}x)e^{-2\pi il\xi}
\ee
for a.e. $(x,\,\xi)\in {\bb R_+}\times \bb R$.
\ed

\br\label{rtheta}
Observe that, given $f\in L^{2}(\bb R_{+})$,
$$\int_{a^j}^{a^{j+1}}\sum\limits_{l\in\bb Z}a^{l}|f(a^{l}x)|^{2}dx=\|f\|_{L^{2}(\bb R_{+})}^{2}<\infty{\mbox{ for }}j\in\bb Z.$$
This implies that  $\sum\limits_{l\in\bb Z}a^{l}|f(a^{l}\cdot)|^{2}<\infty$ a.e. on $\bb R_+$ by the arbitrariness of $j$. Therefore, (\ref{2001}) is well-defined.
\er
\bl\label{onb}~~Let $a$ be as in the general setup. For $m$, $j\in\bb Z$, define $\Lambda_{m}$ as in (\ref{1111}), and $e_{m,j}$ by
$$e_{m,j}(x,\,\xi)=\Lambda_{m}(x)e^{2\pi ij\xi}\mbox{ for }(x,\,\xi)\in \bb R_+\times \bb R.$$
Then

{\rm{(}}i{\rm{)}} $\{\Lambda_m:\,m\in\bb Z\}$ and $\{e_{m,j}:\,m,\,j\in\bb Z\}$ are orthonormal bases for $L^{2}([1,\,a))$ and $L^{2}([1,\,a)\times [0,\,1))$, respectively;

{\rm{(}}ii{\rm{)}} Given $f\in L^2(\bb R_+)$, we have
$$\Theta_{a}f(a^{j}x,\,\xi+m)=e^{2\pi ij\xi}a^{-\frac{j}{2}}\Theta_{a}f(x,\,\xi)$$
for  $j$, $m\in\bb Z$ and a.e. $(x,\,\xi)\in \bb R_+\times \bb R$;

{\rm{(}}iii{\rm{)}} For $j,\,m\in\bb Z$, $f\in L^{2}(\bb R_{+})$,
$$\Theta_{a}(\Lambda_{m}D_{a^{j}}f)(x,\,\xi)=e_{m,j}(x,\,\xi)\Theta_{a}f(x,\,\xi)~{\mbox{ for a.e. }}(x,\,\xi)\in \bb R_+\times \bb R;$$

{\rm{(}}iv{\rm{)}} $\Theta_{a}$-transform is a unitary operator from $ L^{2}(\bb R_{+})$ onto $ L^{2}([1,\,a)\times[0,\,1))$;

{\rm{(}}v{\rm{)}}
\be\label{2222}
\int_{[1,\,a)\times [0,\,1)}|f(x,\,\xi)|^{2}dxd\xi=\sum\limits_{m,j\in \bb Z}\left|\int_{[1,\,a)\times [0,\,1)}f(x,\,\xi)\overline{e_{m,j}(x,\,\xi)}dxd\xi\right|^{2}
\ee for $f\in L^{1}([1,\,a)\times [0,\,1))$.
\el
\bpr By a standard argument, we have {\rm{(}}i{\rm{)}}-{\rm{(}}iii{\rm{)}}.
Next we prove {\rm{(}}iv{\rm{)}} and {\rm{(}}v{\rm{)}}.

{\rm{(}}iv{\rm{)}} It is easy to check that $\Theta_{a}$-transform is a linear and bijective mapping from $L^2(\bb R_+)$ onto $L^2([1,\,a)\times [0,\,1))$. We only need to prove that it is  norm-preserving. For $f\in L^{2}(\bb R_{+})$,
\begin{align*}\left\|\Theta_{a}f\right\|_{L^{2}((1,a)\times(0,1))}^{2}&
=\left\|\sum\limits_{l\in\bb Z}a^{\frac{l}{2}}f(a^{l}x)e^{-2\pi il\xi}\right\|_{L^{2}((1,a)\times(0,1))}^{2} \\&=\int_{1}^{a}dx\int_{0}^{1}\left|\sum\limits_{l\in\bb Z}a^{\frac{l}{2}}f(a^{l}x)e^{-2\pi il\xi}\right|^{2}d\xi\\&
=\int_{1}^{a}\sum\limits_{l\in\bb Z}a^{l}|f(a^{l}x)|^{2}dx\\&
=\|f\|_{L^{2}(\bb R_{+})}^{2}. \end{align*}
This implies that $\Theta_{a}$-transform is norm-preserving.

{\rm{(}}v{\rm{)}} By {\rm{(}}i{\rm{)}}, (\ref{2222}) holds if $f\in L^{2}([1,\,a)\times [0,\,1))$. When $f\in L^{1}([1,\,a)\times [0,\,1))\setminus L^{2}([1,\,a)\times [0,\,1))$, the left-hand side of (\ref{2222}) is infinity. Now we prove by contradiction that the right-hand side of (\ref{2222}) is also infinity. Suppose it is finite, then the function
$$g=\sum\limits_{m,j\in \bb Z}\left(\int_{[1,\,a)\times [0,\,1)}f(x,\,\xi)\overline{e_{m,j}(x,\,\xi)}dxd\xi\right)e_{m,j}$$
belongs to $L^{2}([1,\,a)\times [0,\,1))$ by {\rm{(}}i{\rm{)}}, and thus to $L^{1}([1,\,a)\times [0,\,1))$. It has the same Fourier coefficients as $f$. So $f=g$ by the uniqueness of Fourier coefficients, and thus $f\in L^{2}([1,\,a)\times [0,\,1))$. It is a contradiction.
The proof is completed.
\epr
\br\label{deter} We call the property (ii) the quasi-periodicity of $\Theta_{a}$-transform. By (iv),
 an arbitrary function $F\in L^{2}([1,\,a)\times[0,\,1))$ determines
 a unique  $f\in L^{2}(\bb R_{+})$ in the following way.
Observe that there exists a  unique $\left\{\,c_{m,j}\,\right\}_{m,j\in\bb Z}\in l^2(\bb Z^2)$ such that
$${F(x,\,\xi)=\sum\limits_{m,j\in\bb Z}c_{m,j}e_{m,j}(x,\,\xi)}=\sum\limits_{j\in\bb Z}\left(\sum\limits_{m\in\bb Z}c_{m,\,j}
\Lambda_{m}(x)\right)e^{2\pi ij\xi}$$ for a.e. $(x,\,\xi)\in [1,\,a)\times [0,\,1)$ by (i).
Define $f$ on $\bb R_+$ by {$$f(a^jx)=a^{-\frac j2}\sum\limits_{m\in\bb Z}c_{m,\,-j}
\Lambda_{m}(x){\mbox{ for }}j\in\bb Z{\mbox{ and a.e. }}x\in [1,\,a).$$} Then $$\Theta_{a}f(x,\,\xi)=F(x,\,\xi){\mbox{ for a.e. }}(x,\,\xi)\in [1,\,a)\times [0,\,1).$$
Therefore, we can define $L^2(\bb R_+)$-functions in  $\Theta_{a}$-transform domain.
\er
\bl\label{onbe} Let $a$ and $\Psi$ be as in the general setup. Then
$$
\sum\limits_{l=1}^{L}\sum\limits_{m,j\in \bb Z}\left|\langle f,\,\Lambda_{m}D_{a^{j}}\psi_l\rangle_{L^{2}(\bb R_{+})}\right|^{2}=
\int_{[1,\,a)\times[0,\,1)}\left(\sum\limits_{l=1}^{L}|\Theta_{a}\psi_l(x,\xi)|^{2}\right)
|\Theta_{a}f(x,\xi)|^{2}dxd\xi{\mbox{ for }}f\in L^{2}(\bb R_{+}).
$$
\el
\bpr
Fix $f\in L^{2}(\bb R_{+})$. By  Lemma \ref{onb} (iii) and (iv),  we have
\begin{align}\label{5111}
\sum\limits_{l=1}^{L}
\sum\limits_{m,j\in \bb Z}\left|\langle f,\,\Lambda_{m}D_{a^{j}}\psi_l\rangle_{L^{2}(\bb R_{+})}\right|^{2}
&=\sum\limits_{l=1}^{L}
\sum\limits_{m,j\in \bb Z}
\left|\langle \Theta_{a}f,\,\Theta_{a}\Lambda_{m}D_{a^{j}}
\psi_l\rangle_{L^{2}([1,\,a)\times[0,\,1))}\right|^{2} \non\\
&=\sum\limits_{l=1}^{L}\sum\limits_{m,j\in \bb Z}\left|\langle \Theta_{a}f,\,e_{m,j}\Theta_{a}\psi_l\rangle_{L^{2}
([1,\,a)\times[0,\,1))}\right|^{2}
\non\\
&=\sum\limits_{l=1}^{L}\sum\limits_{m,j\in \bb Z}\left|\int_{[1,\,a)\times[0,\,1)}\overline{\Theta_{a}\psi_l(x,\,\xi)}
\Theta_{a}f(x,\,\xi)\overline{e_{m,j}(x,\,\xi)}dxd\xi\right|^2.
\non
\end{align}
Again applying Lemma \ref{onb} (v) to $\overline{\Theta_{a}\psi_l(x,\,\xi)}
\Theta_{a}f(x,\,\xi)$
leads to
\begin{align*}
  \sum\limits_{l=1}^{L}\sum\limits_{m,j\in \bb Z}\left|\langle f,\,\Lambda_{m}D_{a^{j}}\psi_l\rangle_{L^{2}(\bb R_{+})}\right|^{2}& =\sum\limits_{l=1}^{L}\int_{[1,\,a)\times[0,\,1)}
 \left| \overline{\Theta_{a}\psi_l(x,\,\xi)}
\Theta_{a}f(x,\,\xi)\right|^2
dxd\xi\\
&=\int_{[1,\,a)\times[0,\,1)}\left(\sum\limits_{l=1}^{L}|\Theta_{a}\psi_l(x,\xi)|^{2}\right)
|\Theta_{a}f(x,\xi)|^{2}dxd\xi.
\end{align*}
This finishes the proof.
\epr

\bt\label{complete1} Let $a$ and $\Psi$ be as in the general setup. Then ${\cal MD}(\Psi,\,a)$ is complete in $L^{2}(\bb R_{+})$ if and only if
\be \label{2005}
\sum\limits_{l=1}^{L}\left|\Theta_{a}\psi_{l}(x,\,\xi)\right|^{2}\neq0{\mbox{ for a.e. }}(x,\,\xi)\in [1,\,a)\times[0,\,1).
\ee\et
\bpr
By Lemma \ref{onbe}, for $f\in L^2(\bb R_+)$,
\be \label{2006}
\sum\limits_{l=1}^{L}\sum\limits_{m,j\in \bb Z}\left|\langle f,\,\Lambda_{m}D_{a^{j}}\psi_l\rangle_{L^{2}(\bb R_{+})}\right|^{2}=0~{\mbox{ a.e. on }}\bb R_{+}\ee
if and only if
\be \label{2007}
\left(\sum\limits_{l=1}^{L}|\Theta_{a}\psi_l(x,\xi)|^{2}\right)
|\Theta_{a}f(x,\xi)|^{2}=0~{\mbox{ for a.e. }}~(x,\,\xi)\in [1,\,a)\times[0,\,1).\ee
Observe that ${\cal MD}(\Psi,\,a)$ is complete in $L^{2}(\bb R_{+})$ if and only if $f=0$ is a unique solution to (\ref{2006}) in $L^{2}(\bb R_{+})$. It follows that the completeness of ${\cal MD}(\Psi,\,a)$ in $L^{2}(\bb R_{+})$
is equivalent to $f=0$ being a unique solution to (\ref{2007}) in $L^{2}(\bb R_{+})$. This is in turn  equivalent to the fact that $\Theta_{a}f=0$ is a unique solution to (\ref{2007}) in $L^{2}([1,\,a)\times[0,\,1))$ by Lemma \ref{onb} (iv), which is equivalent to (\ref{2005}). The proof is completed.
\epr

\bt \label{bessel} Let $a$ and $\Psi$ be as in the general setup. Then ${\cal MD}(\Psi,\,a)$ is a Bessel sequence in $L^{2}(\bb R_{+})$ with the Bessel bound $B$ if and only if
\be \label{5112}
\sum\limits_{l=1}^{L}|\Theta_{a}\psi_{l}(x,\xi)|^{2}\leq B ~{\mbox{ for a.e. }}~(x,\,\xi)\in [1,\,a)\times[0,\,1).\ee
\et

\bpr~~Applying Lemma \ref{onbe}, we have
\be \label{5113}\sum\limits_{l=1}^{L}\sum\limits_{m,j\in \bb Z}\left|\langle f,\,\Lambda_{m}D_{a^{j}}\psi_{l}\rangle_{L^{2}(\bb R_{+})}\right|^{2}=
\int_{[1,\,a)\times[0,\,1)}\left(\sum\limits_{l=1}^{L}|\Theta_{a}\psi_{l}
(x,\,\xi)|^2\right)
\left|\Theta_{a}f(x,\,\xi)\right|^2dxd\xi {\mbox{ for }}f\in L^2(\bb R_+).
\ee
So (\ref{5112}) implies that
\be \label{5114}
\sum\limits_{l=1}^{L}\sum\limits_{m,j\in \bb Z}\left|\langle f,\,\Lambda_{m}D_{a^{j}}\psi_{l}\rangle_{L^{2}(\bb R_{+})}\right|^{2}\leq B\int_{[1,\,a)\times[0,\,1)}
|\Theta_{a}f(x,\,\xi)|^{2}dxd\xi
=B\|f\|_{L^{2}_{(\bb R_{+})}}
\ee
for $f\in  L^2(\bb R_+)$ by Lemma \ref{onb} (iv). Thus ${\cal MD}(\Psi,\,a)$ is a Bessel sequence in $L^{2}(\bb R_{+})$ with the Bessel bound $B$.

Now we prove the converse implication by contradiction. Suppose
 ${\cal MD}(\Psi,\,a)$ is a Bessel sequence in $L^{2}(\bb R_{+})$ with the Bessel bound $B$, and $\sum\limits_{l=1}^{L}|\Theta_{a}\psi(\cdot,\cdot)|^{2}> B $ on some $E\subset [1,\,a)\times[0,\,1)$ with $|E|>0$. Take $f$ by
$$\Theta_{a}f(\cdot,\,\cdot)=\chi_{_E}(\cdot,\,\cdot)\mbox{ on }[1,\,a)\times[0,\,1)$$ in (\ref{5113}),
where  $\chi_{_E}$ denotes the characteristic function of $E$.
Then $f$ is well-defined,
$$\|f\|_{L^{2}(\bb R_{+})}^{2}=\int_{[1,\,a)\times[0,\,1)}|\Theta_{a}f(x,\,\xi)|^{2}dxd\xi=|E|$$
by Lemma \ref{onb} (iv), and
\begin{align*}
\sum\limits_{l=1}^{L}\sum\limits_{m,j\in \bb Z}\left|\langle f,\,\Lambda_{m}D_{a^{j}}\psi_{l}\rangle_{L^{2}(\bb R_{+})}\right|^{2}> B|E|
=B\|f\|_{L^{2}_{(\bb R_{+})}}.\end{align*}
It contradicts the fact that ${\cal MD}(\Psi,\,a)$ is a Bessel sequence in $L^{2}(\bb R_{+})$ with the Bessel bound $B$. The proof is completed.
\epr
By an argument similar to  Theorem \ref{bessel}, we have

\bt \label{frame1} Let $a$ and $\Psi$ be as in the general setup.  Then ${\cal MD}(\Psi,\,a)$ is a frame in $L^{2}(\bb R_{+})$ with frame bounds $A$ and $B$ if and only if
$A\leq\sum\limits_{l=1}^{L}|\Theta_{a}\psi_{l}(x,\xi)|^{2}\leq B $ for a.e. $(x,\,\xi)\in [1,a)\times[0,1)$.
\et

\section{$\Theta_{a}$-transform domain expression of duals}
\setcounter{equation}{0}

In this section, we  characterize and express ${\cal MD}$-duals of a general frame ${\cal MD}(\Psi,\,a)$ for $L^{2}(\bb R_{+})$.
And we also study the redundancy of a general frame ${\cal MD}(\Psi,\,a)$ for $L^{2}(\bb R_{+})$. Interestingly, we prove that an arbitrary frame ${\cal MD}(\Psi,\,a)$ for $L^{2}(\bb R_{+})$ is always nonredundant if $L=1$, and is
always redundant if $L>1$ (see Theorem \ref{5-8-5} below).

For convenience, we write
\be\label{D}\mathfrak{D}=\{f\in L^{2}(\bb R_{+}):~\Theta_{a}f\in L^{\infty}([1,\,a)\times[0,\,1))\}.\ee
Then $\mathfrak{D}$ is dense in $L^2(\bb R_+)$
by following Lemma \ref{onb} (iv) and the fact that  $L^{\infty}([1,\,a)\times[0,\,1))$ is dense in $L^{2}([1,\,a)\times[0,\,1))$. This facts will be frequently used in what follows.

Let $a$ and $\Psi$ be as in the general setup, and
 ${\cal MD}(\Psi,\,a)$ be a Bessel sequence in $L^{2}(\bb R_{+})$. We denote by $S$ its frame operator, i.e.,
$$Sf=\sum\limits_{l=1}^{L}\sum\limits_{m, j\in \bb Z}
\langle f,\,\Lambda_{m}D_{a^{j}}\psi_{l}\rangle_{L^{2}(\bb R_{+})}\Lambda_{m}D_{a^{j}}\psi_{l}~~\mbox{ for }~f\in L^{2}(\bb R_{+}).$$
By a standard argument, we have the following lemma which shows that $S$ commutes with the modulation and dilation operators.

\bl\label{5-8-2} Let $a$ and $\Psi$ be as in the general setup. Assume that
 ${\cal MD}(\Psi,\,a)$ is a Bessel sequence in $L^{2}(\bb R_{+})$, and that
 $S$ is its frame operator.  Then
$$S\Lambda_m f=\Lambda_m Sf,~~SD_{a^j}f=D_{a^j}Sf,$$
and thus $S\Lambda_m D_{a^j}f=\Lambda_m D_{a^j}Sf$ for $f\in L^{2}(\bb R_{+})$ and $m,\,j\in\bb Z$.
\el
\bl\label{3-9-1} Let $a$ and $\Psi$ be as in the general setup,  and $\Phi=\{\varphi_1,\,\varphi_2,\,\cdot\cdot\cdot,\,\varphi_L\}\subset L^2(\bb R_+)$. Then
\be\label{391}\sum\limits_{l=1}^{L}\sum\limits_{m,j\in \bb Z}\langle f,\,\Lambda_{m}D_{a^{j}}\psi_l\rangle_{L^{2}(\bb R_{+})}\langle \Lambda_{m}D_{a^{j}}\varphi_l,\,g\rangle_{L^{2}(\bb R_{+})}=\int_{[1,\,a)\times[0,\,1)}\Omega(x,\,\xi)
\Theta_{a}f(x,\,\xi)\overline{\Theta_{a}g(x,\,\xi)}dxd\xi\ee
for $f,\,g\in\mathfrak{D}$, where
$$
\Omega(x,\,\xi)=\sum\limits_{l=1}^{L}\Theta_{a}\varphi_
l(x,\,\xi)\overline{\Theta_{a}\psi_l(x,\,\xi)}.
$$
\el
\bpr Fix $f,\,g\in\mathfrak{D}$. Then $$
\sum\limits_{l=1}^{L}\sum\limits_{m,j\in \bb Z}|\langle f,\,\Lambda_mD_{a^j}\psi_l\rangle|^2<\infty, {\mbox{ and }}\sum\limits_{l=1}^{L}\sum\limits_{m,j\in \bb Z}|\langle g,\,\Lambda_mD_{a^j}\varphi_l\rangle|^2<\infty$$
by Lemma \ref{onbe}, and thus
 the series $$\sum\limits_{l=1}^{L}\sum\limits_{m,j\in \bb Z}\langle f,\,\Lambda_{m}D_{a^{j}}\psi_l\rangle_{L^{2}(\bb R_{+})}\langle \Lambda_{m}D_{a^{j}}\varphi_l,\,g\rangle_{L^{2}(\bb R_{+})}$$ converges absolutely and  is well-defined. By Lemma \ref{onb} (i), (iii) and (iv), we see that
\begin{align*}&\sum\limits_{l=1}^{L}\sum\limits_{m,j\in \bb Z}\langle f,\,\Lambda_{m}D_{a^{j}}\psi_l\rangle_{L^{2}(\bb R_{+})}\langle \Lambda_{m}D_{a^{j}}\varphi_l,\,g\rangle_{L^{2}(\bb R_{+})}\\&=\sum\limits_{l=1}^{L}\sum\limits_{m,j\in \bb Z}\langle \Theta_{a}f,\,\Theta_{a}\Lambda_{m}D_{a^{j}}\psi_{l}\rangle_{L^{2}([1,\,a)\times[0,\,1))}
\langle \Theta_{a}\Lambda_{m}D_{a^{j}}\varphi_{l},\,\Theta_{a}g\rangle_{L^{2}([1,\,a)\times[0,\,1))}\\&=
\sum\limits_{l=1}^{L}\sum\limits_{m,j\in \bb Z}\langle \overline{\Theta_{a}\psi_{l}}\Theta_{a}f,\,e_{m,j}\rangle_{L^{2}([1,\,a)\times[0,\,1))}
\langle e_{m,j},\,\overline{\Theta_{a}\varphi_{l}}\Theta_{a}g\rangle_{L^{2}([1,\,a)\times[0,\,1))}\\&=
\sum\limits_{l=1}^{L}\langle\Theta_{a}f\overline{\Theta_{a}\psi_{l}},\,
\Theta_{a}g\overline{\Theta_{a}\varphi_{l}} \rangle_{L^{2}([1,\,a)\times[0,\,1))}\\&=
\int_{[1,\,a)\times[0,\,1)}
\Omega(x,\,\xi)
\Theta_{a}f(x,\,\xi)\overline{\Theta_{a}g(x,\,\xi)}dxd\xi.
\end{align*}
The proof is completed.\epr
\bl\label{5-8-1} Let $a$ and $\Psi$ be as in the general setup. Assume that
 ${\cal MD}(\Psi,\,a)$ is a Bessel sequence in $L^{2}(\bb R_{+})$, and that
 $S$ is its frame operator. Then,  for $f\in L^{2}(\bb R_{+})$,
\be \label{581} \Theta_{a}Sf(\cdot,\,\cdot)=\left(\sum\limits_{l=1}^{L}\left|\Theta_{a}\psi_{l}
(\cdot,\,\cdot)\right|^{2}\right)\Theta_{a}f(\cdot,\,\cdot)
\ee
a.e. on $[1,\,a)\times[0,\,1)$.
\el
\bpr By Lemma \ref{3-9-1}, we have
$$\langle Sf,\, g\rangle_{L^{2}(\bb R_{+})}=\int_{[1,\,a)\times[0,\,1)}
\left(\sum\limits_{l=1}^{L}\left|\Theta_{a}\psi_{l}
(x,\,\xi)\right|^{2}\right)\Theta_{a}f(x,\,\xi)\overline{\Theta_{a}g(x,\,\xi)}dxd\xi$$
for $f,\,g\in \mathfrak{D}$. Since $\mathfrak{D}$ is dense in $L^{2}(\bb R_{+})$ and ${\cal MD}(\Psi,\,a)$ is a Bessel sequence, by Theorem \ref{bessel} and a  standard argument, it follows that
$$\langle Sf,\, g\rangle_{L^{2}(\bb R_{+})}=\left\langle\left(\sum\limits_{l=1}^{L}\left|\Theta_{a}\psi_{l}
(x,\,\xi)\right|^{2}\right)\Theta_{a}f,\, \Theta_{a}g\right\rangle_{L^{2}([1,\,a)\times[0,\,1))}$$
for $f,\,g\in L^{2}(\bb R_{+})$. Also observing that $\langle Sf,\, g\rangle_{L^{2}(\bb R_{+})}=\left\langle\Theta_{a}Sf,\, \Theta_{a}g\right\rangle_{L^{2}([1,\,a)\times[0,\,1))}$ by Lemma \ref{onb} (iv), we see that
$$\left\langle\Theta_{a}Sf,\, \Theta_{a}g\right\rangle_{L^{2}([1,\,a)\times[0,\,1))}=\left\langle\left(\sum\limits_{l=1}^{L}\left|\Theta_{a}\psi_{l}
\right|^{2}\right)\Theta_{a}f,\, \Theta_{a}g\right\rangle_{L^{2}([1,\,a)\times[0,\,1))}$$
for $f,\,g\in L^{2}(\bb R_{+})$. This implies (\ref{581}) by Lemma \ref{onb} (iv). The proof is completed.
\epr
\bl\label{5-8-4} Let $a$ and $\Psi$ be as in the general setup. Then there exists no Riesz sequence ${\cal MD}(\Psi,\,a)$ in $L^{2}(\bb R_{+})$ whenever $L>1$.
\el
\bpr~By contradiction. Suppose $L>1$ and
${\cal MD}(\Psi,\,a)$ is a Riesz sequence in $L^{2}(\bb R_{+})$. Let $S$ be its frame operator. Then it commutes $\Lambda_m D_{a^j}$ for all $m,\,j\in \bb Z$ by Lemma \ref{5-8-2}. Since $S$ is self-adjoint, invertible and bounded, it follows that
$$S^{-\frac{1}{2}}\Lambda_m D_{a^j}\psi_l=\Lambda_m D_{a^j}S^{-\frac{1}{2}}\psi_l\mbox{ for }~m,\,j\in \bb Z {\mbox{ and }}1\le l\le L.$$
Therefore, ${\cal MD}(S^{-\frac{1}{2}}(\Psi),\,a)$ is an orthonormal  system in $L^{2}(\bb R_{+})$. Write $S^{-\frac{1}{2}}\psi_l=\varphi_l$ for $1\leq l\leq L$. Then
$$\langle\Lambda_{m_1} D_{a^{j_1}}\varphi_{l_1},\, \Lambda_{m_2} D_{a^{j_2}}\varphi_{l_2}\rangle_{L^{2}(\bb R_{+})}=\delta_{m_1,m_2}\delta_{j_1,j_2}\delta_{l_1,l_2}$$
for $m_1,\,m_2,\,j_1,\,j_2\in\bb Z$ and $1\leq l_1,l_2\leq L$,
where the  Kronecker delta is defined by $ \delta_{n,\,m}  = \left\{
           \begin{array}{ll}
            1 &\textrm{if $n=m$}; \\
             0& \textrm{if $n\neq m$}.
           \end{array} \right.$
By Lemma \ref{onb} (iii) and (iv), it is equivalent to
$$\langle e_{m_1,j_1}\Theta_a\varphi_{l_1},\, e_{m_2,j_2}\Theta_a\varphi_{l_2}\rangle_{L^{2}([1,\,a)\times[0,\,1))}=\delta_{m_1,m_2}\delta_{j_1,j_2}\delta_{l_1,l_2}$$
for $m_1,\,m_2,\,j_1,\,j_2\in\bb Z$ and $1\leq l_1,l_2\leq L$, equivalently,
$$\frac{1}{\sqrt{a-1}}\int_{[1,\,a)\times[0,\,1)}\Theta_a\varphi_{l_1}(x,\,\xi)
\overline{\Theta_a\varphi_{l_2}(x,\,\xi)e_{m,j}(x,\,\xi)}dxd\xi=\delta_{m,0}
\delta_{j,0}\delta_{l_1,l_2}$$
for $m,\,j\in\bb Z$ and $1\leq l_1,l_2\leq L$. This is in turn equivalent to
$$\Theta_a\varphi_{l_1}(\cdot,\,\cdot)
\overline{\Theta_a\varphi_{l_2}(\cdot,\,\cdot)}=\delta_{l_1,l_2}~~\mbox{ a.e. on }~[1,\,a)\times[0,\,1)$$
for $1\leq l_1,l_2\leq L$ by the uniqueness of Fourier coefficients. In particular, it implies that
$$|\Theta_a\varphi_{1}(\cdot,\,\cdot)|=|\Theta_a\varphi_{2}(\cdot,\,\cdot)|=1$$
and
$$\Theta_a\varphi_{1}(\cdot,\,\cdot)\overline{\Theta_a\varphi_{2}(\cdot,\,\cdot)}=0$$
a.e. on $[1,\,a)\times[0,\,1)$. This is a contradiction. The proof is completed.
\epr

The following lemma is borrowed from [\ref{lz16},\,Corollary 3.1].
\bl \label{poo} Let $a$ and $\Psi$ be as in the general setup, and $L=1$. Then ${\cal MD}(\Psi,\,a)$ is a Parseval frame for $L^{2}(\bb R_{+})$
if and only if it is an orthonormal basis for $L^{2}(\bb R_{+})$.
\el
\bt \label{dual} Let $a$ and $\Psi$ be as in the general setup, and $\Phi=\{\varphi_1,\,\varphi_2,\,\cdots,\,\varphi_L\}\subset L^2(\bb R_+)$. Assume that ${\cal MD}(\Psi,\,a)$ and ${\cal MD}(\Phi,\,a)$ are Bessel sequences in $L^{2}(\bb R_{+})$. Then
${\cal MD}(\Psi,\,a)$ and ${\cal MD}(\Phi,\,a)$ form a pair of dual frames for $L^{2}(\bb R_{+})$ if and only if
\be \label{5115}\sum\limits_{l=1}^{L}\Theta_{a}\varphi_{l}(x,\,\xi)\overline{\Theta_{a}\psi_{l}(x,\,\xi)}=1{\mbox{ for a.e. }}(x,\,\xi)\in [1,\,a)\times[0,\,1).\ee
\et

\bpr~~Since ${\cal MD}(\Psi,\,a)$ and ${\cal MD}(\Phi,\,a)$ are Bessel sequences in $L^{2}(\bb R_{+})$, and $\mathfrak{D}$ is dense in $L^{2}(\bb R_{+})$, we see that ${\cal MD}(\Psi,\,a)$ and ${\cal MD}(\Phi,\,a)$ form a pair of dual frames for $L^{2}(\bb R_{+})$ if and only if
\be \label{5116}\sum\limits_{l=1}^{L}\sum\limits_{m,j\in \bb Z}\langle f,\,\Lambda_{m}D_{a^{j}}\psi_l\rangle_{L^{2}(\bb R_{+})}\langle \Lambda_{m}D_{a^{j}}\varphi_l,\,g\rangle_{L^{2}(\bb R_{+})}=\langle f,\,g\rangle_{L^{2}(\bb R_{+})}\ee
for $f,\,g\in\mathfrak{D}$. By Lemma \ref{3-9-1} and Lemma \ref{onb} (iv),
(\ref{5116}) is equivalent to
\be \label{5117}\int_{[1,\,a)\times[0,\,1)}\left(\sum\limits_{l=1}^{L}\Theta_{a}\varphi_l(x,\,\xi)
\overline{\Theta_{a}\psi_l(x,\,\xi)}\right)
\Theta_{a}f(x,\,\xi)\overline{\Theta_{a}g(x,\,\xi)}dxd\xi=\int_{[1,\,a)\times[0,\,1)}
\Theta_{a}f(x,\,\xi)\overline{\Theta_{a}g(x,\,\xi)}dxd\xi
\ee
for $f,\,g\in\mathfrak{D}$. Obviously, (\ref{5115}) implies (\ref{5117}). Now we prove the converse implication to finish the proof. Suppose (\ref{5117}) holds. By Theorem \ref{bessel} and the cauchy-Schwarz inequality, we have $\sum\limits_{l=1}^{L}\Theta_{a}\varphi_l\overline{\Theta_{a}\psi_l}\in L^{\infty}([1,\,a)\times(0,\,1))$. It implies that almost every point in $(1,\,a)\times(0,\,1)$ is a Lebesgue point of $\sum\limits_{l=1}^{L}\Theta_{a}\varphi_l\overline{\Theta_{a}\psi_l}$. Arbitrarily fix such a point $(x_0,\,\xi_0)\in (1,\,a)\times(0,\,1)$, and take $f,\,g\in\mathfrak{D}$ in (\ref{5117}) such that
$$\Theta_{a}f=\Theta_{a}g=
\frac{1}{\sqrt{|B((x_{0},\,\xi_{0}),\,\varepsilon)|}}
\chi_{_B((x_{0},\,\xi_{0}),\,\varepsilon)}.$$
on $[1,\,a)\times[0,\,1)$ with  $B((x_{0},\,\xi_{0}),\,\varepsilon)\subset(1,\,a)\times(0,\,1)$ and $\varepsilon>0$, where $B((x_{0},\,\xi_{0}),\,\varepsilon)$ denotes the $\varepsilon$-neighborhood of $(x_{0},\,\xi_{0})$. Then $f$ and $g$ are
well-defined by Lemma \ref{onb} (iv), and we obtain that
\be \label{5118}\frac{1}{|B((x_{0},\,\xi_{0}),\,\varepsilon)|}\int_{B((x_{0},\,\xi_{0}),\,\varepsilon)}
\sum\limits_{l=1}^{L}\Theta_{a}\varphi_{l}(x,\,\xi)
\overline{\Theta_{a}\psi_{l}(x,\,\xi)}dxd\xi=1.
\ee
Letting $\varepsilon\rightarrow 0$ in (\ref{5118}) leads to
$$\sum\limits_{l=1}^{L}\Theta_{a}\varphi_{l}(x_{0},\,\xi_{0})
\overline{\Theta_{a}\psi_{l}(x_{0},\,\xi_{0})}=1.$$
This implies (\ref{5115}) by the arbitrariness of $(x_{0},\,\xi_{0})$.
The proof is completed.
\epr

Now we turn to the expression of ${\cal MD}$-duals. Let $a$ and $\Psi$ be as in the general setup, ${\cal MD}(\Psi,\,a)$ be a frame for $L^{2}(\bb R_{+})$ and $S$ be its frame operator. By Lemma \ref{5-8-2}, $S\Lambda_m D_{a^j}=\Lambda_m D_{a^j}S$, and thus $S^{-1}\Lambda_m D_{a^j}=\Lambda_m D_{a^j}S^{-1}$ for $m,\,j\in\bb Z$. So ${\cal MD}(\Psi,\,a)$ and its canonical dual $S^{-1}({\cal MD}(\Psi,\,a))$ share the same dilation-and-modulation structure, that is,
$$S^{-1}({\cal MD}(\Psi,\,a))={\cal MD}(S^{-1}(\Psi),\,a).$$ The following theorem gives its canonical dual window and all ${\cal MD}$-dual windows in $\Theta_a$ transform domain.

\bt\label{5-8-3} Let $a$ and $\Psi$ be as in the general setup, and ${\cal MD}(\Psi,\,a)$ be a frame for $L^{2}(\bb R_{+})$. Then

{\rm{(}}i{\rm{)}}~its canonical dual ${\cal MD}(S^{-1}(\Psi),\,a)$ is given by
$$\Theta_a S^{-1}\psi_{l}(\cdot,\,\cdot)=\frac{\Theta_a\psi_{l}(\cdot,\,\cdot)}
{\sum\limits_{l=1}^{L}\left|\Theta_{a}\psi_{l}
(\cdot,\,\cdot)\right|^{2}}~~\mbox{ a.e. on }~[1,\,a)\times[0,\,1){\mbox{ for }}1\le l\le L;$$

{\rm{(}}ii{\rm{)}}~a dilation-and-modulation system ${\cal MD}(\Phi,\,a)$ with
$\Phi=\{\varphi_1,\,\varphi_2,\,\cdots,\,\varphi_L\}$ is a dual frame of ${\cal MD}(\Psi,\,a)$ if and only if $\Phi$ is defined by
\be \label{582}\Theta_a \varphi_{l}(\cdot,\,\cdot)=\frac{\Theta_a\psi_{l}(\cdot,\,\cdot)
\left(1-\sum\limits_{l=1}^{L}\overline{\Theta_a\psi_{l}
(\cdot,\,\cdot)}X_{l}(\cdot,\,\cdot)\right)}
{\sum\limits_{l=1}^{L}\left|\Theta_{a}\psi_{l}
(\cdot,\,\cdot)\right|^{2}}+X_{l}(\cdot,\,\cdot)~~\mbox{ a.e. on }~[1,\,a)\times[0,\,1),\ee
where $X_{l}\in L^{\infty}([1,\,a)\times[0,\,1))$ with $1\leq l\leq L$.
\et
\bpr {\rm{(}}i{\rm{)}} Since $S$ is an invertible and bounded operator on $L^{2}(\bb R_{+})$, we have
\be \label{583}
\Theta_a f(\cdot,\,\cdot)=\left(\sum\limits_{l=1}^{L}\left|\Theta_{a}\psi_{l}
(\cdot,\,\cdot)\right|^{2}\right)\Theta_{a}S^{-1}f(\cdot,\,\cdot)~~\mbox{ for }~f\in L^{2}(\bb R_{+})
\ee
by Lemma \ref{5-8-1}. Replacing $f$ by $\psi_{l}$ in (\ref{583}) with $1\leq l\leq L$, we have {\rm{(}}i{\rm{)}}.

{\rm{(}}ii{\rm{)}} Sufficiency. Suppose $\Phi$ is given by (\ref{582}). Then ${\cal MD}(\Phi,\,a)$ is a Bessel sequence in $L^{2}(\bb R_{+})$ by Theorem \ref{bessel}. By a simple computation, we have
$$\sum\limits_{l=1}^{L}\Theta_a \varphi_{l}(\cdot,\,\cdot)\overline{\Theta_a \psi_{l}(\cdot,\,\cdot)}=1~~\mbox{ a.e. on }~[1,\,a)\times[0,\,1).$$
It follows that ${\cal MD}(\Phi,\,a)$ is a dual frame of ${\cal MD}(\Psi,\,a)$ by Theorem \ref{dual}.

Necessity. Suppose ${\cal MD}(\Phi,\,a)$ is a dual frame of ${\cal MD}(\Psi,\,a)$. Then
$$\sum\limits_{l=1}^{L}\overline{\Theta_a \psi_{l}(\cdot,\,\cdot)}\Theta_a \varphi_{l}(\cdot,\,\cdot)=1~~\mbox{ a.e. on }~[1,\,a)\times[0,\,1)$$
by Theorem \ref{dual}, and $\Theta_a \varphi_{l}\in L^{\infty}([1,\,a)\times[0,\,1))$. So we have (\ref{582}) with $X_{l}=\Theta_a \varphi_{l}, 1\leq i\leq L$. The proof is completed.
\epr

The following theorem shows that the cardinality $L$ of $\Psi$ determines whether or not a frame ${\cal MD}(\Psi,\,a)$ is redundant. If $L=1$, there exists no redundant frame ${\cal MD}(\Psi,\,a)$ for $L^{2}(\bb R_{+})$. If $L>1$, there exists no nonredundant frame ${\cal MD}(\Psi,\,a)$ for $L^{2}(\bb R_{+})$.

\bt\label{5-8-5} Let $a$ and $\Psi$ be as in the general setup, and ${\cal MD}(\Psi,\,a)$ be a frame for $L^{2}(\bb R_{+})$. Then ${\cal MD}(\Psi,\,a)$ is a Riesz basis for $L^{2}(\bb R_{+})$ if and only if $L=1$.
\et
\bpr The necessity is an immediate consequence of Lemma \ref{5-8-4}. Now we show the sufficiency. Suppose $L=1$. From the proof of Lemma \ref{5-8-4},
$$S^{-\frac{1}{2}}({\cal MD}(\Psi,\,a))={\cal MD}(S^{-\frac{1}{2}}(\Psi),\,a).$$
So ${\cal MD}(S^{-\frac{1}{2}}(\Psi),\,a)$ is a Parseval frame for $L^{2}(\bb R_{+})$ since ${\cal MD}(\Psi,\,a)$ is a frame for $L^{2}(\bb R_{+})$. It leads to that ${\cal MD}(S^{-\frac{1}{2}}(\Psi),\,a)$ is an orthonormal basis by Lemma \ref{poo}. This is equivalent to the fact that ${\cal MD}(\Psi,\,a)$ is a Riesz basis for $L^{2}(\bb R_{+})$. The proof is completed.
\epr
\section{Some examples}
\setcounter{equation}{0}
Theorems \ref{bessel}, \ref{frame1} and \ref{5-8-3} provide us with an easy method to construct ${\cal MD}$-dual frame pairs for $L^{2}(\bb R_{+})$. This section focus on some examples. They show that we can construct ${\cal MD}$-dual frame pairs for $L^{2}(\bb R_{+})$ with good properties such as dual windows having bounded supports and certain smoothness.

\begin{example}\label{611-1} Let $c$ be a finitely supported  sequence defined on $\bb Z$, and its Fourier transform
$$\hat{c}(\xi)=\sum\limits_{l\in\bb Z}c_{l}e^{-2\pi il\xi}$$
have no zero on $[0,\,1)$. Define $\psi\in L^{2}(\bb R_{+})$ by
$$\Theta_a\psi(x,\,\xi)=\hat{c}(\xi)~~~~\mbox{ for }~(x,\,\xi)\in[1,\,a)\times[0,\,1).$$
Then $\psi$ is a step function and of bounded support by the definition of $\Theta_a$, and ${\cal MD}(\psi,\,a)$ is a frame for $L^{2}(\bb R_{+})$ by Theorem \ref{frame1} since $|\hat{c}(\xi)|$ have positive lower and upper bounds due to its continuity and having no zeros on $[0,\,1)$. It follows that ${\cal MD}(\psi,\,a)$ has a unique ${\cal MD}$-dual window $S^{-1}\psi$ defined by
$$\Theta_aS^{-1}\psi(x,\,\xi)=\frac{1}{\sum\limits_{l\in\bb Z}\overline{c_{l}}e^{2\pi il\xi}}~~~~\mbox{ for }~(x,\,\xi)\in[1,\,a)\times[0,\,1)$$
by Theorems \ref{5-8-3} and \ref{5-8-5}. Observe that, if at least two $c_{l}$ are nonzero, we have
$$\frac{1}{\sum\limits_{l\in\bb Z}\overline{c_{l}}e^{2\pi il\xi}}=\sum\limits_{l\in\bb Z}d_{l}e^{-2\pi il\xi}$$
with $d$ being infinitely supported. It follows that the dual window $S^{-1}\psi$  is of unbounded support by the definition of $\Theta_a$, although $\psi$ is of bounded support.
\end{example}

The following  example shows that it is possible for us to obtain multi-window ${\cal MD}$-dual frame pairs for $L^{2}(\bb R_{+})$ with each window being of bounded support.

\begin{example}\label{611-2} Let $L>1$, $m_1,\,m_2,\,\cdot\cdot\cdot,\,m_L$ be trigonometric polynomials satisfying
$$|m_1(\xi)|^2+|m_2(\xi)|^2+\cdot\cdot\cdot+|m_L(\xi)|^2=1~~~\mbox{ for }\xi\in [0,\,1).$$
Define $\Psi=\{\psi_1,\,\psi_2,\,\cdots,\,\psi_L\}$ by
$$\Theta_a\psi_l(x,\,\xi)=m_l(\xi)~~~~\mbox{ for }~(x,\,\xi)\in[1,\,a)\times[0,\,1).$$
Then ${\cal MD}(\Psi,\,a)$ is a frame for $L^{2}(\bb R_{+})$ and every $\psi_l$ is of bounded support by an argument similar to Example \ref{611-1}.
Define $\Phi=\{\varphi_1,\,\varphi_2,\,\cdots,\,\varphi_L\}$ by
\be \label{6111}\Theta_a\phi(x,\,\xi)=m_l(\xi)\left(1-\sum\limits_{l=1}^{L}\overline{m_l(\xi)}
X_l(x,\,\xi)\right)+X_l(x,\,\xi)~~~~\mbox{ for a.e. }~(x,\,\xi)\in[1,\,a)\times[0,\,1)\ee
with $X_l\in L^{\infty}([1,\,a)\times[0,\,1))$. Then ${\cal MD}(\Psi,\,a)$ and ${\cal MD}(\Phi,\,a)$ form a pair of dual frames for $L^{2}(\bb R_{+})$ by Theorem \ref{5-8-3}. If, in addition, we require
\be \label{6112}X_l(x,\,\xi)=\sum\limits_{j\in\bb Z}d_{l,j}(x)e^{-2\pi ij\xi}~~~~\mbox{ for a.e. }~(x,\,\xi)\in[1,\,a)\times[0,\,1).\ee
with $\left\{d_{l,j}(\cdot)\right\}_{j\in\bb Z}$ is a finitely supported sequence of functions on $[1,\,a)$ for each $1\leq l\leq L$, then each $\varphi_l$ with $1\leq l\leq L$ is of bounded support by (\ref{6111}) and the definition of $\Theta_a$.
\end{example}

\begin{example}\label{611-3} Let $L\geq 1$, $\Psi=\{\psi_1,\,\psi_2,\,\cdots,\,\psi_L\}$ be a finite subset of  $L^{2}(\bb R_{+})$, and $supp(\psi_l)\subset[1,\,a)$.
Assume that
$$\sum\limits_{l=1}^{L}|\psi_l(x)|^2=1~~~~\mbox{ for a.e. }~x\in[1,\,a).$$
Define $\Phi=\{\varphi_1,\,\varphi_2,\,\cdots,\,\varphi_L\}$ by
\be \label{6113}\Theta_a\varphi_l(x,\,\xi)=\psi_l(x)\left(1-\sum\limits_{l=1}^
{L}\overline{\psi_l(x)}
X_l(x,\,\xi)\right)+X_l(x,\,\xi)~~~~\mbox{ for a.e. }~(x,\,\xi)\in[1,\,a)\times[0,\,1)\ee
with $X_l\in L^{\infty}([1,\,a)\times[0,\,1))$. Then ${\cal MD}(\Psi,\,a)$ and ${\cal MD}(\Phi,\,a)$ form a pair of dual frames for $L^{2}(\bb R_{+})$ by Theorem \ref{5-8-3}. In particular, if $X_l, 1\leq l\leq L$ are required as in (\ref{6112}), in addition, each $\varphi_l$ with $1\leq l\leq L$ is of bounded support.
\end{example}

In Examples \ref{611-2} and \ref{611-3}, $\Theta_a\psi_l, 1\leq l\leq L$, are
defined by univariate functions. Next we give a relatively more general example.

\begin{example} Assume that $c_0(x)$ and $c_1(x)$ are two real-valued measurable functions defined on $[1,\,a]$, and that there exist two positive constants $A$ and $B$ such that
$$A\leq |c_0(x)|+|c_1(x)|\leq B~~~~\mbox{ for  }~x\in[1,\,a].$$
Define $\Psi=\{\psi_1,\,\psi_2\}\subset L^{2}(\bb R_{+})$ by
$$\Theta_a\psi_1(x,\,\xi)=c_0(x)+c_1(x)e^{-4\pi i\xi}$$
$$ \Theta_a\psi_2(x,\,\xi)  = \left\{
           \begin{array}{ll}
            2i\sqrt{c_0(x)c_1(x)}\sin2\pi\xi &\textrm{if $c_0(x)c_1(x)\geq 0$}; \\
            \\
             2\sqrt{-c_0(x)c_1(x)}\cos2\pi\xi& \textrm{if $c_0(x)c_1(x)<0$}
           \end{array} \right.$$
for a.e. $(x,\,\xi)\in[1,\,a]\times[0,\,1).$  Then
$$ \psi_1(x)  = \left\{
           \begin{array}{ll}
            c_0(x) &\textrm{if $1\leq x\leq a$}; \\
            \\
            a^{-1}c_1(a^{-2}x) &\textrm{if $a^{2}\leq x\leq a^{3}$};
            \\
            \\
             0& \textrm{otherwise},
           \end{array} \right.$$
and
$$ \psi_2(x)  = \left\{
           \begin{array}{ll}
            a^{\frac{1}{2}}\sqrt{c_0(ax)c_1(ax)} &\textrm{if $a^{-1}\leq x\leq 1$}; \\
            \\
            -a^{-\frac{1}{2}}\sqrt{c_0(a^{-1}x)c_1(a^{-1}x)} &\textrm{if $a\leq x\leq a^{2}$ and $c_0(a^{-1}x)c_1(a^{-1}x)\geq 0$};
            \\
            \\
            a^{-\frac{1}{2}}\sqrt{-c_0(a^{-1}x)c_1(a^{-1}x)} &\textrm{if $a\leq x\leq a^{2}$ and $c_0(a^{-1}x)c_1(a^{-1}x)< 0$};\\
            \\
             0& \textrm{otherwise},
           \end{array} \right.$$
 and  $$|\Theta_a\psi_1(x,\,\xi)|^2
+|\Theta_a\psi_2(x,\,\xi)|^2=(|c_0(x)|+|c_1(x)|)^2$$
by a simple computation and the definition of $\Theta_a$.
It follows that
\be\label{6114} A^2\leq |\Theta_a\psi_1(x,\,\xi)|^2
+|\Theta_a\psi_2(x,\,\xi)|^2\leq B^2\ee
for a.e. $(x,\,\xi)\in[1,\,a]\times[0,\,1)$, and thus ${\cal MD}(\Psi,\,a)$
is a frame for $L^{2}(\bb R_{+})$ by Theorem \ref{frame1}.
Obviously, $\psi_1$ and  $\psi_2$  are real-valued and of bounded support.

Now we check the ${\cal MD}$-duals of ${\cal MD}(\Psi,\,a)$. Define $\Phi=\{\varphi_1,\,\varphi_2\}$ by
$$\Theta_a\varphi_l(x,\,\xi)=\frac{\Theta_a\psi_l(x,\,\xi)
\left(1-\overline{\Theta_a\psi_1(x,\,\xi)}X_1(x,\,\xi)-
\overline{\Theta_a\psi_2(x,\,\xi)}X_2(x,\,\xi)\right)}
{\left(\left|c_0(x)\right|+\left|c_1(x)\right|\right)^2}+X_l(x,\,\xi)$$
for $1\leq l\leq 2$ and a.e.$(x,\,\xi)\in[1,\,a]\times[0,\,1)$ with $X_1,\,X_2\in L^{\infty}([1,\,a]\times[0,\,1))$.  Then ${\cal MD}(\Psi,\,a)$ and ${\cal MD}(\Phi,\,a)$ form a pair of dual frames for $L^{2}(\bb R_{+})$
by Theorem \ref{5-8-3}. If
$$X_l(x,\,\xi)=\sum\limits_{j\in\bb Z}d_{l,j}(x)e^{-2\pi ij\xi}~~~~\mbox{ for a.e. }~(x,\,\xi)\in[1,\,a]\times[0,\,1).$$
with $\left\{d_{l,j}(\cdot)\right\}_{j\in\bb Z}$ being a finitely supported sequence of real-valued function on $[1,\,a]$ for each $1\leq l\leq 2$, then $\varphi_1$ and $\varphi_2$
are also real-valued and of bounded support. Also
we can  obtain $\Phi$ with good smoothness by choosing good
$X_1$ and $X_2$. For example, if we make further assumption that $c_0(x)$, $c_1(x)$  and $\sqrt{|c_0(x)c_1(x)|}$ are $k$-th continuously differentiable on $(1,\,a)$, that
$c_0(1)c_1(1)=c_0(a)c_1(a)=0$, and $c_0(x)c_1(x)>0$ for $x\in(1,\,a)$. Then
$\psi_1$ and $\psi_2$ are continuous functions on $\bb R_{+}$ and $k$-th continuously differentiable on $(1,\,a)\cup (a^2,\,a^3)$ and $(a^{-1},\,1)\cup (a,\,a^2)$,
respectively. In this case,  if we further require that $|c_0(x)|+|c_1(x)|$ is a constant  on $[1,\,a]$, and $X_1(x,\,\xi)$ and $X_2(x,\,\xi)$ satisfy
$$X_1(x,\,\xi)=\sum\limits_{j\in\bb Z}d_{1,j}e^{-2\pi ij\xi}$$
and
$$X_2(x,\,\xi)=\sum\limits_{j\in\bb Z}d_{2,j}e^{-2\pi ij\xi}$$
for $\xi\in [0,\,1)$ with $\left\{d_{1,j}\right\}_{j\in\bb Z}$ and $\left\{d_{2,j}\right\}_{j\in\bb Z}$ being two finitely supported real number sequences. Then $\varphi_1$ and $\varphi_2$ are real-valued, of bounded support, and have the same continuity
and differentiability as $\psi_1$ and $\psi_2$.
\end{example}

\end{document}